\documentclass[11pt,a4paper]{article}
\usepackage{amsfonts}
\usepackage{amsmath, amssymb, graphics, setspace}
\usepackage{graphicx}
\usepackage{url}
\usepackage{listings}
\usepackage{indentfirst}

\newtheorem{teor}{Theorem}[section]
\newcommand{\mathsym}[1]{{}}
\newcommand{\unicode}[1]{{}}

\title{On the numerical Picard iterations method with collocations for the IVP}
\author{E. Scheiber\thanks{email: scheiber@unitbv.ro}}
\date{}

\begin{document}
\maketitle
\begin{abstract}
Some variants of the numerical Picard iterations method are presented to solve an IVP
for an ordinary differential system. The term numerical
emphasizes that a numerical solution is computed. The method consists in replacing the
right hand side of the differential system by  Lagrange interpolation polynomials
followed by successive approximations. In the case when the number of interpolation point
is fixed a convergence result is given. Finally some numerical experiments are reported.
\end{abstract}



\section{Introduction}
The paper presents variants of the Picard iterations to solve an initial value problem (IVP)
for ordinary differential equations. On subintervals the right hand side of the differential system
is replaced by  Lagrange interpolation polynomials  on subintervals and then successive approximations are used.
The interpolation nodes are the images of a set of reference points.
The number of these reference points can be fixed or variable, i.e. increasing number \cite{2}.

When the number of reference nodes is fixed the approximations of the solution of the IVP
are computed by collocations. A convergence result is given. This case appears in \cite{10}, p. 211.
In \cite{1} the spectral deferred correction is defined adding a term to the iteration formula and 
the convergence of that method is proved.

If the number of reference points increases then  the values of the unknown function are 
determined iteratively \cite{2}.

We use the terminology numerical Picard iterations to emphasize that the
method builds a numerical solution. For an IVP the usual Picard iterations are exemplified
with Computer Algebra code in \cite{7}.

For stiff problems the Picard iterations method is treated in \cite{8}, \cite{9}.

There is another approach to the numerical Picard iterations for an IVP, where the approximations 
are a linear form of Chebyshev polynomials \cite{5}, \cite{6}, \cite{1}.

In the last section some results of our computational experiences are presented.

\section{Numerical Picard iterations}

Let the IVP be
\begin{eqnarray}
\dot{y}(x) &=& f(x,y(x)), \qquad x\in[x_0,x_f],\label{sac1}\\
y(x_0) &=& y^0, \label{sac2}
\end{eqnarray}
where the function $f:[x_0,x_f]\times\mathbb{R}^N\rightarrow\mathbb{R}^N$ has the components
$f=(f^1,\ldots,f^N)$. 

In $\mathbb{R}^N$ for $y=(y^1,\ldots,y^N)$ we shall use the norm $\|y\|=\max_{1\le j\le N}|y^j|.$ 

We assume that $f$ is Lipschitz continuous, i.e. there exists 
$L>0$ such that
$$
|f^{\mu}(x,y_1)-f^{\mu}(x,y_2)|\le L\sum_{j=1}^N|y_1^j-y_2^j|\qquad \forall\ y_1,y_2\in\mathbb{R}^N,\ \mu\in\{1,2,\ldots,N\}
$$
and consequently
$$
\|f(x,y_1)-f(x,y_2)\|\le \tilde{L}\|y_1-y_2\|,
$$
where $\tilde{L}=NL.$

The IVP (\ref{sac1})-(\ref{sac2}) may be reformulated as the integral equation
\begin{equation}\label{sac3}
y(x)=y^0+\int_{x_0}^xf(s,y(s))\mathrm{d}s.
\end{equation}

Within these hypotheses the problem (\ref{sac1})-(\ref{sac2}) or (\ref{sac3}) has a unique solution. This solution
may be obtained with the Picard iterations 
\begin{eqnarray*}
y^{(n+1)}(x)&=&y^0+\int_{x_0}^xf(s,y^{(n)}(s))\mathrm{d}s,\quad n\in\mathbb{N},\\
y^{(0)}(x) &=& y^0,
\end{eqnarray*}
for $x\in[x_0,x_f].$ The sequence $(y^{(n)}(x))_{k\in\mathbb{N}}$ converges uniformly in $[x_0,x_f]$
to the solution of IVP.

Let  $M\in\mathbb{N}^*,\ h=\frac{x_f-x_0}{M}$ and the mesh be defined as
$x_i=x_0+ih,\ i\in\{0,1,\ldots,M\}.$
The numerical solution is given by the sequence  $u_h=(u_0,u_1,\ldots,u_M),$ where each $u_i=u(x_i)$
is an approximation of $y(x_i).$

If $u_i$ was computed, on the interval $[x_i,x_{i+1}]$ the function  $f(s,y(s))$
under the integral in
\begin{equation}\label{sac14}
y(x)=y(x_i)+\int_{x_i}^xf(s,y(s))\mathrm{d}s
\end{equation}
will be replaced by a Lagrange interpolation polynomial
\begin{equation}\label{sac4}
u(x)=u(x_i)+\int_{x_i}^xL(\mathbb{P}_{m-1};x_{i,1},x_{i,2},\ldots,x_{i,m};f(\cdot,u(\cdot)))(s)\mathrm{d}s,\quad x\in[x_i,x_{i+1}].
\end{equation}
The interpolation nodes  $x_i\le x_{i,1}<x_{i,2}<\ldots<x_{i,m}\le x_{i+1}$ are fixed by a certain rule.
The used notation states the interpolation constraints
$$
L(\mathbb{P}_{m-1};x_{i,1},x_{i,2},\ldots,x_{i,m};f(\cdot,u(\cdot)))(x_{i,j})=f(x_{i,j},u(x_{i,j})),\qquad j\in\{1,2,\ldots,m\}.
$$
From (\ref{sac4}) we deduce
\begin{equation}\label{sac5}
u(x)=u(x_i)+\sum_{j=1}^m\left(\int_{x_i}^xl_j(s)\mathrm{d}s\right)f(x_{i,j},u(x_{i,j})),
\end{equation}
where $(l_j)_{1\le j\le m}$ are the Lagrange fundamental polynomials
\begin{equation}\label{sac6}
l_j(x)=\frac{(x-x_{i,1})\ldots(x-x_{i,j-1})(x-x_{i,j+1})\ldots(x-x_{i,m})}
{(x_{i,j}-x_{i,1})\ldots(x_{i,j}-x_{i,j-1})(x_{i,j}-x_{i,j+1})\ldots(x_{i,j}-x_{i,m})}.
\end{equation}

\section{Picard iterations with a fixed reference set}

In (\ref{sac5}) the values
$$
u(x_{i,1}), u(x_{i,2}), \ldots, u(x_{i,m})
$$
are unknown. To compute these vectors the collocation method will be used.

Choosing $x:=x_{i,k}$ in (\ref{sac5}) we get
\begin{equation}\label{sac7}
u(x_{i,k})=u(x_i)+\sum_{j=1}^m\left(\int_{x_i}^{x_{i,k}}l_j(s)\mathrm{d}s\right)f(x_{i,j},u(x_{i,j})),
\quad k\in\{1,2,\ldots,m\}.
\end{equation}
The relations (\ref{sac7}) form a nonlinear system with the unknowns
 $u(x_{i,1}),\ldots,u(x_{i,m}) \in \underbrace{\mathbb{R}^N\times\ldots\times\mathbb{R}^N}_m \simeq \mathbb{R}^{mN}.$

In order to simplify and provides a unitary approach to the computation of the integrals from (\ref{sac7})
we fix the nodes $\xi_1<\xi_2<\ldots<\xi_m$ within an interval $[a,b].$
We call these nodes the reference interpolation nodes.
The function 
$$
\varphi_i(\xi)=x_i+\frac{h}{b-a}(\xi-a)
$$ 
maps the interval $[a,b]$ into $[x_i,x_{i+1}].$ For any $i\in\{0,1,\ldots,M-1\}$ the nodes 
$x_{i,j}$ will be defined as
$$
x_{i,j}=\varphi_i(\xi_j),\qquad \forall\ j\in\{1,2,\ldots,m\}.
$$
If $s=\varphi_i(\xi)$ then
$$
l_j(s)=\frac{(\xi-\xi_1)\ldots(\xi-\xi_{j-1})(\xi-\xi_{j+1})\ldots(\xi-\xi_m)}
{(\xi_j-\xi_1)\ldots(\xi_j-\xi_{j-1})(\xi_j-\xi_{j+1})\ldots(\xi_j-\xi_m)}=\tilde{l}_j(\xi)
$$
and
$$
\int_{x_i}^{x_{i,k}}l_j(s)\mathrm{d}s=\frac{h}{b-a}\int_a^{\xi_k}\tilde{l}_j(\xi)\mathrm{d}\xi.
$$
Denoting $w_{j,k}=\frac{1}{b-a}\int_a^{\xi_k}\tilde{l}_j(\xi)\mathrm{d}\xi$
the nonlinear system (\ref{sac7}) becomes
\begin{equation}\label{sac11}
u(x_{i,k})=u(x_i)+h\sum_{j=1}^mw_{j,k}f(x_{i,j},u(x_{i,j})),
\quad k\in\{1,2,\ldots,m\}.
\end{equation}

In order to prove the existence of a solution of the nonlinear system we shall use
a simplified notation
$u(x_{i,k})=u_k,\ k\in\{1,2,\ldots,m\}.$ Then the system (\ref{sac11}) is written as
\begin{equation}\label{sac12}
u_k=u(x_i)+h\sum_{j=1}^mw_{j,k}f(x_{i,j},u_j),
\quad k\in\{1,2,\ldots,m\}.
\end{equation}
The operator
$$
\Phi=(\Phi_k)_{1\le k\le m}\quad\mbox{where}\quad
\Phi_k:\underbrace{\mathbb{R}^N\times\ldots\times\mathbb{R}^N}_m\rightarrow\mathbb{R}^N
$$
is defined by
$$
\Phi_k(u)=u(x_i)+h\sum_{j=1}^mw_{j,k}f(x_{i,j},u_j),\qquad u=(u_1,\ldots,u_m).
$$
The used norm in $\underbrace{\mathbb{R}^N\times\ldots\times\mathbb{R}^N}_m$ will be
$$
\|u\|=\|(u_1,\ldots,u_m)\|=\sum_{j=1}^m\|u_j\|.
$$

If $u=(u_1,\ldots,u_m)$ and $v=(v_1,\ldots,v_m)$ then following equality is valid
$$
\Phi_k(u)-\Phi_k(v)=h\sum_{j=1}^mw_{j,k}\left(f(x_{i,j},u_j)-f(x_{i,j},v_j)\right),\quad k\in\{1,2,\ldots,m\}.
$$
Then
$$
\|\Phi_k(u)-\Phi_k(v)\|\le h \tilde{L}\sum_{j=1}^m|w_{j,k}|\|u_j-v_j\|
$$
and
$$
\sum_{k=1}^m\|\Phi_k(u)-\Phi_k(v)\|\le h \tilde{L}\sum_{k,j=1}^m|w_{j,k}|\|u_j-v_j\|.
$$
If $\omega=\max_{1\le j\le m}\sum_{k=1}^m|w_{j,k}|$ the above inequality gives
$$
\|\Phi(u)-\Phi(v)\|\le h\omega \tilde{L}\|u-v\|.
$$
Following theorem is a consequence of the above:
\begin{teor}
For $h$ small enough ($h<\frac{1}{\omega L}$) the nonlinear system (\ref{sac11}) has a unique solution.
\end{teor}

In the hypothesis of the above theorem, the nonlinear system (\ref{sac11}) may be solved using the 
successive approximation method
\begin{eqnarray}
u^{(n+1)}(x_{i,k}) &=& u(x_i)+h\sum_{j=1}^mw_{j,k}f(x_{i,j},u^{(n)}(x_{i,j})),
\quad n\in\mathbb{N}\label{sac8}\\
u^{(0)}(x_{i,k}) &=& u(x_i) \label{sac9}
\end{eqnarray}
for $k\in\{1,2,\ldots,m\}.$ The sequences
$$
u^{(n)}_{i,j}\stackrel{\mathrm{def}}{=}u^{(n)}(x_{i,j}),\qquad n\in\mathbb{N},\quad j\in\{1,2,\ldots,m\}
$$
will converge to the solution of the system (\ref{sac11}).

The iterative relations (\ref{sac8}) can be written in matrix form
\begin{equation}\label{sac18}
[u^{(n+1)}_{i,1}\ u^{(n+1)}_{i,2}\ldots u^{(n+1)}_{i,m}]=
\end{equation}
$$
=\underbrace{[u_i\ u_i\ldots u_i]}_{m}
+h[f_{i,1}\ f_{i,2}\ldots f_{i,m}]\left(\begin{array}{ccc}
w_{1,1} &\ldots & w_{1,m}\\
w_{2,1}  &\ldots & w_{2,m}\\
\vdots & \ddots & \vdots\\
w_{m,1} &\ldots & w_{m,m}
\end{array}\right),
$$
where $f_{i,j}=f(x_{i,j},u^{(n)}_{i,j}), \ j\in\{1,\ldots,m\}.$

Denoting $u^{(n)}_i=(u^{(n)}_{i,j})_{1\le j\le m}$ the iterations stop when the following 
condition is fulfilled $\|u^{(n)}_i-u^{(n-1)}_i\|< \varepsilon,$
where $\varepsilon>0$ is a tolerance.
The initial approximations are chosen as $u^{(0)}_{i,j}=u(x_i)$ for any $j\in\{1,2,\ldots,m\}.$

This method to solve the nonlinear system (\ref{sac11}) leads to an approximation
to the solution of the IVP in the most right node which may differ from  $x_{i+1}$.
We point out two variants of the computations:
\begin{itemize}
\item
We change the initial mesh such that $x_{i+1}$ will be the most right node
($x_{i+1}=\varphi_i(\xi_m)$)
and the computation continue in the interval $[x_{i+1},x_{i+1}+h].$
In this case we have
$$
u_{i+1}\stackrel{\mathrm{def}}{=}u(x_{i+1})=u^{(n)}_{i,m}.
$$
\item
In (\ref{sac4}) we set $x:=x_{i+1}=\varphi_i(b)$ and
$$
u_{i+1}\stackrel{\mathrm{def}}{=}u(x_{i+1})=u(x_i)+\frac{h}{b-a}\sum_{j=1}^m\left(\int_a^b\tilde{l}_j(\xi)\mathrm{d}\xi\right)f(x_{i,j},u^{(n)}_{i,j}).
$$
In this way $m$ new integrals must be computed additionally.
\end{itemize}
With the new notations we have $u_0\stackrel{\mathrm{def}}{=}u(x_0)=y^0.$

\vspace*{0.3cm}
The coefficients $w_{j,k}$ do not depend on the computation interval.
We highlight some cases when these coefficients may be easily computed.
 
\subsection*{Some particular cases}
\begin{enumerate}
\item
Equidistant nodes. If 
$\xi_j=\frac{j-1}{m-1},\ j\in\{1,2,\ldots,m\},$
then
$$
w_{j,k}=\int_0^{\xi_k}\tilde{l}_j(\xi)\mathrm{d}\xi=\frac{(-1)^{m-j}}{(j-1)!(m-j)!}
\int_0^{\frac{k-1}{m-1}}\prod_{\stackrel{\mu=0}{\mu\not=j-1}}^{m-1}\left((m-1)\xi-\mu\right)\mathrm{d}\xi.
$$
The following Mathematica code computes these coefficients:

\begin{doublespace}
\noindent\(\pmb{\text{Wcoeff}[\text{j$\_$},\text{k$\_$},\text{m$\_$}]\text{:=}}\\
\pmb{\text{Module}[\{x,w\},}\\
\pmb{w=\text{Integrate}[\text{Product}[\text{If}[i\neq j-1,(m-1)x-i,1],\{i,0,m-1\}],}\\
\pmb{\{x,0,(k-1)/(m-1)\}];(-1){}^{\wedge}(m-j) w/((j-1)! (m-j)!)]}\)
\end{doublespace}

The results obtained for $m=2$ are

\begin{doublespace}
\noindent\(\pmb{\text{MatrixForm}[\text{Table}[\text{Wcoeff}[j,k,2],\{k,1,2\},\{j,1,2\}]]}\)
\end{doublespace}

\begin{doublespace}
\noindent\(\left(
\begin{array}{cc}
 0 & 0 \\
 \frac{1}{2} & \frac{1}{2} \\
\end{array}
\right)\)
\end{doublespace}

Because $x_{i,1}=x_i$ \c{s}i $x_{i,2}=x_{i+1}$ the recurrence formula (\ref{sac8})-(\ref{sac9})  becomes
\begin{eqnarray*}
u^{(n+1)}_{i+1}&=& u_i+\frac{h}{2}\left(f(x_i,u_i)+f(x_{i+1},u^{(n)}_{i+1})\right);\\
u^{(0)}_{i+1} &=& u_i.
\end{eqnarray*}

For $m=3$ the results are

\begin{doublespace}
\noindent\(\pmb{\text{MatrixForm}[\text{Table}[\text{Wcoeff}[j,k,3],\{k,1,3\},\{j,1,3\}]]}\)
\end{doublespace}

\begin{doublespace}
\noindent\(\left(
\begin{array}{ccc}
 0 & 0 & 0 \\
 \frac{5}{24} & \frac{1}{3} & -\frac{1}{24} \\
 \frac{1}{6} & \frac{2}{3} & \frac{1}{6} \\
\end{array}
\right)\)
\end{doublespace}

In this case $x_{i,1}=x_i, x_{i,2}=\frac{1}{2}(x_i+x_{i+1})\stackrel{\mathrm{def}}{=}x_{i+\frac{1}{2}}, x_{i,3}=x_{i+1}$ 
and the recurrence formulas (\ref{sac8})-(\ref{sac9}) become
\begin{eqnarray*}
u^{(n+1)}_{i+\frac{1}{2}} &=& u_i+
h\left(\frac{5}{24}f(x_i,u_i)+\frac{1}{3}f(x_{i+\frac{1}{2}},u^{(n)}_{i+\frac{1}{2}}))
-\frac{1}{24}f(x_{i+1},u^{(n)}_{i+1})\right);\\
u^{(n+1)}_{i+1} &=& u_i+ 
\frac{h}{6}\left(f(x_i,u_i)+4f(x_{i+\frac{1}{2}},u^{(n)}_{i+\frac{1}{2}})
+f(x_{i+1},u^{(n)}_{i+1})\right);\\
u^{(0)}_{i+\frac{1}{2}} &=& u_i;\\
u^{(0)}_{i+1} &=& u_i.
\end{eqnarray*}
In matrix form the above relations are
$$
\left[\begin{array}{l}
u_i \\
u^{(n+1)}_{i+\frac{1}{2}} \\
u^{(n+1)}_{i+1} 
\end{array}\right]=
\left[\begin{array}{l}
u_i \\
u_i \\
u_i 
\end{array}\right]+h\left(\begin{array}{ccc}
0 & 0 & 0\\
\frac{5}{24} & \frac{1}{3} & -\frac{1}{24} \\
\frac{1}{6} & \frac{2}{3} & \frac{1}{6}
\end{array}\right)
\left[\begin{array}{l}
f(x_i,u_i) \\
f(x_{i+\frac{1}{2}},u^{(n)}_{i+\frac{1}{2}})\\
f(x_{i+1},u^{(n)}_{i+1})
\end{array}\right].
$$ 
Transposing the above equality we get the form corresponding to (\ref{sac18}).

To compute $u_{i+1}$ we observe that for $m=2$ the trapezoidal rule, while for $m=3$
the Simpson integration formula are used.
\item
Chebyshev points of second kind $\xi_j=\cos{\frac{(j-1)\pi}{m-1}},\ j\in\{1,\ldots,m\}.$
Then
$$
w_{j,k}=\frac{1}{2}\int_{-1}^{\xi_k}\tilde{l}_j(\xi)\mathrm{d}\xi=
\frac{(-1)^{j-1}2^{m-3}\gamma_j}{m-1}\int_{-1}^{\xi_k}\prod_{k=1,k\not=j}^m(\xi-\xi_k)\mathrm{d}\xi
$$
with $\gamma_j=\left\{\begin{array}{lcl}
0.5 & \mbox{if} & j\in\{1,m\}\\
1 & \mbox{if} & j\in\{2,\ldots,m-1\}
\end{array}\right..$
\item
The nodes are the roots of an orthogonal polynomial.
Now we suppose that the polynomial $p_m(\xi)=\prod_{j=1}^m(\xi-\xi_j)$ is orthogonal to $\mathbb{P}_{m-1},$
the set of polynomials of  degree at most $m-1,$ with the weight  $\rho(\xi)$ on the interval 
$I=[a,b].$ In this case the Lagrange fundamental polynomials
$\tilde{l}_j(\xi),\ j\in\{1,\ldots,m\}$ are orthogonal.
\begin{itemize}
\item
If $\rho(\xi)=1, I=[a,b]$ then $p_m(\xi)=\frac{m!}{(2m)!}\frac{\mathrm{d}^m}{\mathrm{d}^m\xi}(\xi-a)^m(\xi-b)^m$
is the Laguerre polynomial. For $a=0, b=1$ and $m=1$ following results are obtained
\begin{eqnarray*}
p_1(\xi) &=& \xi-\frac{1}{2} \quad \Rightarrow\quad \xi_1=\frac{1}{2}\\
w_{1,1} &=& \int_0^{\frac{1}{2}}\mathrm{d}\xi=\frac{1}{2}\\
u^{(n+1)}_{i+\frac{1}{2}}  &=& u_i+\frac{h}{2}f(x_{i+\frac{1}{2}},u^{(n)}_{i+\frac{1}{2}})
\end{eqnarray*}

Again we observe that $u(x_{i+\frac{1}{2}})$ is computed using the rectangle rule in the right hand side of
 (\ref{sac4}).
\item
The Chebyshev polynomials  $p_m(\xi)=\frac{1}{2^{m-1}}\cos(m\arccos{\xi}),\ m\in\mathbb{N},$ are orthogonal with the weight
$\rho(\xi)=\frac{1}{\sqrt{1-\xi^2}}$ in $I=[-1,1].$

The nodes will be
$$
\xi_j=\cos{\frac{(2j-1)\pi}{2m}} \quad \Rightarrow\quad x_{i,j}=x_i+\frac{h}{2}(\xi_j+1),\ j\in\{1,2,\ldots,m\}.
$$
The biggest node is $x_{i,1}.$ 
The Lagrange fundamental polynomials are
$$
\tilde{l}_j(\xi)=\frac{2^{m-1}}{m}(-1)^{j-1}\sin{\frac{(2j-1)\pi}{2m}}
\prod_{\stackrel{\mu=1}{\mu\not=j}}^m\left(\xi-\cos{\frac{(2\mu-1)\pi}{2m}}\right)
$$
and
$$
w_{j,k}=\frac{2^{m-2}}{m}(-1)^{j-1}\sin{\frac{(2j-1)\pi}{2m}}
\int_{-1}^{\cos{\frac{(2k-1)\pi}{2m}}}\prod_{\stackrel{\mu=1}{\mu\not=j}}^m\left(\xi-\cos{\frac{(2\mu-1)\pi}{2m}}\right)
\mathrm{d}\xi
$$
The integral can be analytically computed but it involves rounding errors.
\end{itemize}

\end{enumerate}

\subsection{The convergence of the method}

The function $f(x,y(x))$ being continuous there exists a constant $K_1>0$ such that
$$
\max_{1\le\mu\le N}\max_{x\in[x_0,x_f]}|f^{\mu}(x,y(x))|\le K_1.
$$
We also suppose that the function $f(x,y)$ are continuous partial derivatives of order $m$ for any  
$x\in[x_0,x_f].$ There exists $K_m>0$ such that
 $$
\max_{1\le\mu\le N}\max_{x\in[x_0,x_f]}\left|\frac{\mathrm{d}^mf^{\mu}(x,y(x))|}{\mathrm{d}x^m}\right|\le K_m.
$$
In any interval $[x_i,x_{i+1}]$ the following equality is valid
$$
f^{\mu}(x,y(x))-L(\mathbb{P}_{m-1};x_{i,1},x_{i,2},\ldots,x_{i,m};f^{\mu}(\cdot,y(\cdot)))(x)=
$$
$$
=\frac{1}{m!}\prod_{j=1}^m(x-x_{i,j})
\left.\frac{\mathrm{d}^mf^{\mu}(x,y(x))}{\mathrm{d}x^m}\right|_{x=\eta_{\mu}}
$$
where $\eta_{\mu}\in[x_i,x_{i+1}].$ 

\noindent
We denote by  $R^{\mu}(x)$ the right hand side and then
 $\max_{x\in[x_i,x_{i+1}]}|R^{\mu}(x)|\le \frac{K_m}{m!}h^m.$
If $R(x)=(R^1(x),\ldots,R^N(x)) $ then (\ref{sac14})  implies the vectorial relaton
\begin{equation}\label{sac15}
y(x)=y(x_i)+\int_{x_i}^xL(\mathbb{P}_{m-1};x_{i,1},x_{i,2},\ldots,x_{i,m};f^{\mu}(\cdot,y(\cdot)))(s)\mathrm{d}s+
\int_{x_i}^xR(s)\mathrm{d}s
\end{equation}
and $\|\int_{x_i}^xR(s)\mathrm{d}s\|\le\frac{K_m}{m!}h^{m+1}.$

We make the following notations
$$
\begin{array}{lcll}
e_i &=& \|y(x_i)-u_i\|, & i\in\{0,1,\ldots,M\};\\
r^{(n)}_{i,j} &=& \|y(x_{i,j})-u^{(n)}_{i,j}\|,& j\in\{1,2,\ldots,m\};\\
r^{(n)}_i &=& \max_{1\le j\le m}r^{(n)}_{i,j}. & \\
\end{array}
$$
and additionally
$$
w= \max\left\{\max_{1\le j,k\le m}|w_{j,k}|,\max_{1\le j\le m}\frac{1}{b-a}\left|\int_a^b\tilde{l}_j(\xi)\mathrm{d}\xi\right|\right\},
\qquad \tilde{w}=m w.
$$
We emphasize that $n$ represents the number of iterations on an interval $[x_i,x_{i+1}].$
This number differs from one interval to another. For simplicity we omitted the index $i$ when $n$ is written.

Several times the following theorem will be used
\begin{teor}
If  $(z_k)_{k\in\mathbb{N}}$ is a sequence of nonnegative numbers such that
$$
z_{k+1}\le a z_k+b\quad \forall\ k\in\mathbb{N}\quad \mbox{\c{s}i}\quad a,b>0,\ a\not=1,
$$ 
then
$$
z_k\le a^k z_0+ b\frac{a^k-1}{a-1}, \qquad \forall\ k\in\mathbb{N}.
$$
The above inequality implies: if $a>1$ then $z_k\le a^k\left(z_0+\frac{b}{a-1}\right)$ and if
$a<1$ then $z_k\le a^kz_0+\frac{b}{1-a}.$
\end{teor}

In the beginning we determine an evaluation for $r^{(n)}_i.$

For $n=0$ the equalities hold:
$$
y(x_{i,j})-u^{(0)}_{i,j}=y(x_{i,j})-u_i=\left(y(x_{i,j})-y(x_i)\right)+\left(y(x_i)-u_i\right)=
$$
$$
=\int_{x_i}^{x_{i,j}}f(s,y(s))\mathrm{d}s+\left(y(x_i)-u_i\right)
$$
and then we deduce
$$
r^{(0)}_{i,j}\le K_1h +e_i,\ \forall\ j\in\{1,2,\ldots,m\}\quad\Rightarrow\quad r^{(0)}_i\le K_1h +e_i.
$$
If $n>0,$ for $x=x_{i,k}$ the equality (\ref{sac15}) may be written as
\begin{equation}\label{sac16}
y(x_{i,k})=y(x_i)+h\sum_{j=1}^mw_{j,k}f(x_i,j,y(x_{i,j}))+\int_{x_i}^{x_{i,k}}R(s)\mathrm{d}s.
\end{equation}
Subtracting  (\ref{sac8}) from (\ref{sac16})we obtain
$$
y(x_{i,k})-u^{(n+1)}_{i,k}=
$$
$$
=y(x_i)-u_i+h\sum_{j=1}^mw_{j,k}\left(f(x_{i,j},y(x_{i,j}))-f(x_{i,j},u^{(n)}_{i,j})\right)+
\int_{x_i}^{x_{i,k}}R(s)\mathrm{d}s.
$$
It follows that
$$
r^{(n+1)}_{i,k}\le e_i+h\tilde{L}\tilde{w}r^{(n)}_i+\frac{K_m}{m!}h^{m+1}\quad\Rightarrow\quad r^{(n+1)}_i\le
e_i+h\tilde{L}\tilde{w}r^{(n)}_i+\frac{K_m}{m!}h^{m+1}
$$
If $h$ is small enough ($h\tilde{L}\tilde{w}<1$) then
$$
r^{(n)}_i\le (h\tilde{L}\tilde{w})^nr^{(0)}_i+\frac{1}{1-h\tilde{L}\tilde{w}}\left(e_i+\frac{K_m}{m!}h^{m+1}\right)\le
$$
$$
\le(h\tilde{L}\tilde{w})^n(K_1h +e_i)+\frac{1}{1-h\tilde{L}\tilde{w}}\left(e_i+\frac{K_m}{m!}h^{m+1}\right)=
$$
\begin{equation}\label{sac17}
= \left((h\tilde{L}\tilde{w})^n+\frac{1}{1-h\tilde{L}\tilde{w}}\right)e_i+
h^{n+1}(\tilde{L}\tilde{w})^nK_1+\frac{K_mh^{m+1}}{m!(1-h\tilde{L}\tilde{w})}.
\end{equation}

Evaluating $e_i$ we distinguish two cases depending on the definition of  $u_{i+1}:$

$$
u_{i+1}=u^{(n)}_{i,m}=u_i+h\sum_{j=1}^mw_{j,m}f(x_{i,j},u^{(n-1)}_{i,j}),\qquad (x_{i+1}=\varphi_i(\xi_m))
$$
or
$$
u_{i+1}=u_i+\frac{h}{b-a}\sum_{j=1}^m\left(\int_a^b\tilde{l}_j(\xi)\mathrm{d}\xi\right)f(x_{i,j},u^{(n)}_{i,j}),
\qquad (x_{i+1}=\varphi_i(b)).
$$
Corresponding to the two cases, from (\ref{sac15}) we obtain the equalities
$$
y(x_{i+1})=y(x_i)+h\sum_{j=1}^mw_{j,m}f(x_{i,j},y(x_{i,j}))+\int_{x_i}^{x_{i,m}}R(s)\mathrm{d}s
$$
and respectively
$$
y(x_{i+1})=y(x_i)+\frac{h}{b-a}\sum_{j=1}^m\left(\int_a^b\tilde{l}_j(\xi)\mathrm{d}\xi\right)f(x_{i,j},y(x_{i,j}))+
$$
$$
+\int_{x_i}^{x_{i+1}}R(s)\mathrm{d}s.
$$
Computing $y(x_{i+1})-u_{i+1}$ it results
$$
y(x_{i+1})-u_{i+1}=y(x_i)-u_i+h\sum_{j=1}^mw_{j,m}\left(f(x_{i,j},y(x_{i,j}))-u^{(n-1)}_{i,j}\right)+
$$
$$
+\int_{x_i}^{x_{i,m}}R(s)\mathrm{d}s,
$$
respectively
$$
y(x_{i+1})-u_{i+1}=y(x_i)-u_i+\frac{h}{b-a}\sum_{j=1}^m\left(\int_a^b\tilde{l}_j(\xi)\mathrm{d}\xi\right)
\left(f(x_{i,j},y(x_{i,j}))-u^{(n)}_{i,j}\right)+
$$
$$
+\int_{x_i}^{x_{i+1}}R(s)\mathrm{d}s.
$$
It follows that
$$
e_{i+1}\le e_i+h\tilde{L}\tilde{w}r^{(n-1)}_i+\frac{K_m}{m!}h^{m+1}
$$
and
$$
e_{i+1}\le e_i+h\tilde{L}\tilde{w}r^{(n)}_i+\frac{K_m}{m!}h^{m+1}.
$$
We remark that between the two estimates only the upper index of  $r_i$ differs. This justifies
that in the second case $m$ additional integrals must be computed.

From hereon it is sufficient to consider only the first case. Using
(\ref{sac17}) we obtain
$$
e_{i+1}\le e_i+h\tilde{L}\tilde{w}\left(
 \left((h\tilde{L}\tilde{w})^{n-1}+\frac{1}{1-h\tilde{L}\tilde{w}}\right)e_i+
h^{n}(\tilde{L}\tilde{w})^{n-1}K_1+\frac{K_mh^{m+1}}{m!(1-h\tilde{L}\tilde{w})}
\right)+\frac{K_m}{m!}h^{m+1}=
$$
$$
=e_i\left(1+(h\tilde{L}\tilde{w})^n+\frac{h\tilde{L}\tilde{w}}{1-h\tilde{L}\tilde{w}}\right)+
h^{n+1}(\tilde{L}\tilde{w})^nK_1+\frac{K_mh^{m+1}}{m!(1-h\tilde{L}\tilde{w})}.
$$
Because
$h\tilde{L}\tilde{w}< 1 \ \Rightarrow\ (h\tilde{L}\tilde{w})^n\le h\tilde{L}\tilde{w}$
the above inequality becames
$$
e_{i+1}\le e_i\left(1+h\tilde{L}\tilde{w}+\frac{h\tilde{L}\tilde{w}}{1-h\tilde{L}\tilde{w}}\right)+
h^2\tilde{L}\tilde{w}K_1+\frac{K_mh^{m+1}}{m!(1-h\tilde{L}\tilde{w})}.
$$
Consequently
$$
e_i\le \left(1+h\tilde{L}\tilde{w}+\frac{h\tilde{L}\tilde{w}}{1-h\tilde{L}\tilde{w}}\right)^i\left(e_0+
\frac{h^2\tilde{L}\tilde{w}K_1+\frac{K_mh^{m+1}}{m!(1-h\tilde{L}\tilde{w})}}
{h\tilde{L}\tilde{w}+\frac{h\tilde{L}\tilde{w}}{1-h\tilde{L}\tilde{w}}}\right)\le
$$
$$
\le e^{i\left(h\tilde{L}\tilde{w}+\frac{h\tilde{L}\tilde{w}}{1-h\tilde{L}\tilde{w}}\right)}
\left(e_0+
\frac{h\tilde{L}\tilde{w}K_1+\frac{K_mh^m}{m!(1-h\tilde{L}\tilde{w})}}
{\tilde{L}\tilde{w}+\frac{\tilde{L}\tilde{w}}{1-h\tilde{L}\tilde{w}}}\right).
$$

Because $e_0=0,$ from the above inequality it results that: 
$$
\max_{1\le i\le M}e_i \le e^{(x_f-x_0)\tilde{L}\tilde{w}\left(1+\frac{1}{1-h\tilde{L}\tilde{w}}\right)}
\left(\frac{h\tilde{L}\tilde{w}K_1+\frac{K_mh^m}{m!(1-h\tilde{L}\tilde{w})}}
{\tilde{L}\tilde{w}+\frac{\tilde{L}\tilde{w}}{1-h\tilde{L}\tilde{w}}}\right)\rightarrow 0,
$$
for $h\searrow 0\ \Leftrightarrow\ M\rightarrow\infty.$
This proves the convergence of the method.

\section{Picard iterations with a variable reference set}

We shall keep some of the above introduced notations and we shall define those that differ.

Let $a\le \xi_1^m<\xi_2^m<\ldots<\xi_m^m\le b$ be the roots of the polynomial $p_m(x),$
where $(p_m)_{m\in\mathbb{N}}$ is a sequence of orthogonal polynomials with the weight 
$\rho\in L_2[a,b]$ on the interval $[a,b].$ It is assumed that
$\frac{1}{\rho}\in L_2[a,b],$ too. These are requirements of the convergence theorem \cite{2}.

If $\varphi_i$ is the affine function transforming $[a,b]$ onto $[x_i,x_{i+1}]$
then the nodes  are introduced 
$$
x_{i,j}^m=\varphi_i(\xi_j^m),\qquad j\in\{1,2,\ldots,m\},\quad m\in\mathbb{N}^*.
$$
For $x\in[x_i,x_{i+1}],$ we define
$$
u^{m+1}(x)=u_i+\int_{x_i}^xL(\mathbb{P}_{m-1};x_{i,1}^m,x_{i,2}^m,\ldots,x_{i,m}^m;f(\cdot,u^m(\cdot)))(s)\mathrm{d}s=
$$
$$
=u_i+\sum_{j=1}^m\left(\int_{x_i}^xl_j^m(s)\mathrm{d}s\right)f(x_{i,j}^m,u(x_{i,j}^m))=
$$
$$
=u_i+\frac{h}{b-a}\sum_{j=1}^m\left(\int_a^{\zeta}\tilde{l}_j^m(\xi)\mathrm{d}\xi\right)f(x_{i,j}^m,u(x_{i,j}^m)),
$$
where $\zeta=\varphi_i^{-1}(x)$ and
$$
\tilde{l}_j^m(\xi)=\frac{(\xi-\xi_1^m)\ldots(\xi-\xi_{j-1}^m)(\xi-\xi_{j+1}^m)\ldots(\xi-\xi_m^m)}
{(\xi_j^m-\xi_1^m)\ldots(\xi_j^m-\xi_{j-1}^m)(\xi_j^m-\xi_{j+1}^m)\ldots(\xi_j^m-\xi_m^m)}.
$$
The vectors $u_{i,j}^m$ are defined iteratively 
$$
u_{i,1}^1=u_i,
$$
\begin{eqnarray*}
u_{i,1}^2 &=& u_i+\frac{h}{b-a}\left(\int_a^{\xi_1^2}\tilde{l}_1^1(\xi)\mathrm{d}\xi\right) f(x_{i,1}^1,u_{i,1}^1)=\\
 &=& u_i+\frac{h}{b-a}(\xi_1^2-a)f(x_{i,1}^1,u_{i,1}^1);\\
u_{i,2}^2 &=& u_i+\frac{h}{b-a}\left(\int_a^{\xi_2^2}\tilde{l}_1^1(\xi)\mathrm{d}\xi\right) f(x_{i,1}^1,u_{i,1}^1)=\\
 &=& u_i+\frac{h}{b-a}(\xi_2^2-a)f(x_{i,1}^1,u_{i,1}^1).
\end{eqnarray*}
It was taken into account that $\tilde{l}_1(\xi)=1.$ As a rule
$$
u_{i,k}^{m+1}=u^{m+1}(x_{i,k}^{m+1})=u_i+
\frac{h}{b-a}\sum_{j=1}^m\left(\int_a^{x_{i,k}^{m+1}}\tilde{l}_j^m(\xi)\mathrm{d}\xi\right)f(x_{i,j}^m,u(x_{i,j}^m)),
$$
for $k\in\{1,2,\ldots,m+1\}$ \c{s}i $m\in\mathbb{N}^*.$

We must compute
$$
u_{i+1}^{m+1}=u^{m+1}(x_{i+1})=u_i+
\frac{h}{b-a}\sum_{j=1}^m\left(\int_a^b\tilde{l}_j^m(\xi)\mathrm{d}\xi\right)f(x_{i,j}^m,u(x_{i,j}^m)),
$$
too.

\noindent
The computation of the vectors  $u_{i,k}^{m+1},\ k\in\{1,2,\ldots,m+1\},\ u_{i+1}^{m+1}$ can be written in matrix form. 
For simplicity we denote
\begin{eqnarray*}
w_{j,k}&=& \int_a^{x_{i,k}^{m+1}}\tilde{l}_j^m(\xi)\mathrm{d}\xi,\quad j\in\{1,\ldots,m\},\ k\in\{1,\ldots,m+1\},\\
w_j &=& \int_a^b\tilde{l}_j^m(\xi)\mathrm{d}\xi,\quad  j\in\{1,\ldots,m\}
\end{eqnarray*}
and the matrix
$$
W=\frac{h}{b-a}\left(\begin{array}{cccc}
w_{1,1} & w_{2,1} & \ldots & w_{m,1} \\
w_{1,2} & w_{2,2} & \ldots & w_{m,2} \\
\vdots & & & \vdots \\
w_{1,k+1} & w_{2,k+1} & \ldots & w_{m,k+1} \\
w_1 & w_2 & \ldots & w_m
\end{array}\right)\in M_{m+2,m}(\mathbb{R})
$$
$$
F=[f(x_{i,1}^m,u_{i,1}^m),f(x_{i,2},u_{i,2}^m),\ldots,f(x_{i,m}^m,u_{i,m}^m)]\in M_{N,m}(\mathbb{R})
$$
The following equality holds
$$
[u_{i,1}^{m+1},u_{i,2}^{m+1},\ldots,u_{i,m+1}^{m+1},u_{i+1}^{m+1}]^T={\underbrace{[u_i,\ldots,u_i]}_{m+2}}^T+W\cdot F^T.
$$

For an imposed tolerance $\varepsilon>0,$ the iterations occurs until the condition 
 $\|u_{i+1}^{m+1}-u_i^{m}\|<\varepsilon$ is fulfilled.
The initial approximation is  $u_{i+1}^1=u_i.$ 
If the above condition is fulfilled then we set $u_{i+1}=u_{i+1}^{m+1}.$

A convergence result is given in \cite{2}.


\section{Stiff problems}

From (\ref{sac3}), if $s=x_0+h\sigma$ then
$$
y(x)=y(x_0)+h\int_0^{\frac{x-x_0}{h}}f(x_0+h\sigma,y(x_0+h\sigma)\mathrm{d}\sigma,
$$
with $x\in[x_0,x_0+h]\ \Leftrightarrow \ \sigma\in[0,1].$

Setting
$$
y(x_0+h\sigma)=y_0+hv(\sigma)
$$
we derive that $v(0)=0$ and 
$$
\frac{\mathrm{d}v(\sigma)}{\mathrm{d}\sigma}=f(x_0+h\sigma,y_0+hv(\sigma))\quad\Leftrightarrow\quad
v(s)=\int_0^sf(x_0+h\sigma,y_0+hv(\sigma))\mathrm{d}\sigma.
$$
Following \cite{8}, \cite{9}, by the \textit{stabilization principle}, the solution of the partial differential system
\begin{equation}\label{sac100}
\frac{\partial w(s,t)}{\partial t}=-w(s,t)+\int_0^sf(x_0+h\sigma,y_0+h w(\sigma,t))\mathrm{d}\sigma
\end{equation}
has the property, cf. \cite{8}, \cite{9},
\begin{equation}\label{sac103}
\lim_{t\rightarrow\infty}\|w(s,t)-v(s)\|=0,\qquad \mbox{for}\quad s\in[0,1].
\end{equation}

We give a numerical solution to find an approximation of the solution of (\ref{sac100}).

Let be $\tau>0$ and the sequence $t^n=n\tau,\ n\in\mathbb{N}.$  The equation (\ref{sac100}) may
be rewritten as
$$
\frac{\partial e^tw(s,t)}{\partial t}=e^t\int_0^sf(x_0+h\sigma,y_0+h w(\sigma,t))\mathrm{d}\sigma
$$
and integrating from $n\tau$ to $(n+1)\tau$ it results
\begin{equation}\label{sac101}
w(s,t^{n+1})=e^{-\tau} w(s,t^n)+e^{-(n+1)\tau}
\int_{n\tau}^{(n+1)\tau}e^{\eta}\left(\int_0^sf(x_0+h\sigma,y_0+h w(\sigma,\eta))\mathrm{d}\sigma\right)\mathrm{d}\eta.
\end{equation}

Without changing the notation for $w,$ we substitute in (\ref{sac101}) $f(x_0+h\sigma,y_0+h w(\sigma,\eta))$ 
by a Lagrange interpolation polynomial 
\begin{equation}\label{sac102}
w(s,t^{n+1})=e^{-\tau} w(s,t^n)+
\end{equation}
$$
+e^{-(n+1)\tau}
\int_{n\tau}^{(n+1)\tau}e^{\eta}\left(L(\mathbb{P}_{m-1}\xi_1,\ldots,\xi_m;f(x_0+h\ \cdot,y_0+hw(\cdot,\eta))
\mathrm{d}\sigma\right)\mathrm{d}\eta=
$$
$$
=e^{-\tau} w(s,t^n)+e^{-(n+1)\tau}\sum_{j=1}^m
\int_{n\tau}^{(n+1)\tau}e^{\eta}\left(\int_0^s
f(x_0+h\xi_j,y_0+h w(\xi_j,\eta))l_j(\sigma)\mathrm{d}\sigma\right)\mathrm{d}\eta,
$$
where $0=\xi_1<\xi_2<\ldots<\xi_m=1.$

We denote $w^n(s)=w(s,t^n)$  and in the right hand side of (\ref{sac102}) we take $w(\xi_j,\eta)=w^n(\xi_j),$
for any $j\in\{1,2,\ldots,m\}$ and $\eta\in[n\tau,(n+1)\tau].$ Then  
$$
w^{n+1}(s)=e^{-\tau}w^n(s)+(1-e^{-s})\sum_{j=1}^mf(x_0+h\xi_j,y_0+h w^n(\xi_j))\int_0^sl_j(\sigma)\mathrm{d}\sigma.
$$
Denoting $w^n_j=w^n(\xi_j),$ for $s=\xi_k,\ k\in\{1,2,\ldots,m\}$ we obtain the iterative relations
$$
w^{n+1}_k=e^{-\tau}w^n_k+(1-e^{-\tau})\sum_{j=1}^mf(x_0+h\xi_j,y_0+h w^n_j)\int_0^{\xi_k}l_j(\sigma)\mathrm{d}\sigma.
$$
The iterations occurs until the stopping condition $\max_{1\le j\le m}\|w^{n+1}_j-w^n_j\|<\varepsilon$
is fulfilled. Here $\varepsilon>0$ is a tolerance.
According to (\ref{sac103}) we consider $v(1)=w^{n+1}_j$ and the procedure continues with $u_{i+1}=u_i+h w^{n+1}_m.$

\section{Numerical experiences}

Choosing adequate values for $M$, tolerance and
 the maximum allowed iterations number there are obtained acceptable results.

Using computer programs based on these methods we solved the following IVPs:
\begin{enumerate}
\item
 (\cite{3}, p. 234)
$$
\left\{\begin{array}{lcl}
\dot{y} & = & y\frac{4(x+2)^3-y}{(x+2)^4-1},\qquad x\in[0,1],\\
y(0) &=& 15
\end{array}\right.
$$
with the solution $y(x)=1+(x+2)+(x+2)^2+(x+2)^3.$

For $M=5$ and the tolerance $\varepsilon=10^{-5}$ the maximum error $\max_{0\le i\le M}\|y(x_i)-u_i\|$ and
the number of calling the function $f$ are given in  Table \ref{tab1}.

\begin{table} [h]
\begin{center}
\begin{tabular}{|c|c|c|c|}
\hline\hline
\multicolumn{2}{|c|}{Fixed equidistant} & \multicolumn{2}{c|}{Variable reference set} \\
\multicolumn{2}{|c|}{reference set $m=3$} & \multicolumn{2}{c|}{}\\
\hline
 Error & $N_f$ & Error & $N_f$ \\
 \hline\hline
 1.82591e-08 & 75 & 8.94274e-08 & 99 \\
 \hline\hline
\end{tabular}
\end{center}

\vspace*{-0.5cm}
\caption{Results for Example 1.}\label{tab1}
\end{table}

\item
 (\cite{3}, p. 244)
$$
\left\{\begin{array}{lclcl}
\dot{y_1} & = & y_2,&\quad& y_1(0)=1, \qquad x\in[0,x_f],\\
\dot{y_2} & = & -\frac{y_1}{r^3},&\quad& y_2(0)=0, \\
\dot{y_3} & = & y_4,&\quad& y_3(0)=0, \\
\dot{y_4} & = & -\frac{y_3}{r^3},&\quad& y_4(0)=1, 
\end{array}\right.
$$
where $r=\sqrt{y_1^2+y_3^2}$ and with the solution $y_1=\cos{x}, y_2=-\sin{x}, y_3=\sin{x}, y_4=\cos{x}.$

The results of our numerical experiments are listed in Table \ref{tab2}.
\begin{table} [h]
\begin{center}
\begin{tabular}{|c|c|c|c|c|c|c|}
\hline\hline
& & & \multicolumn{2}{|c|}{Fixed equidistant}  & \multicolumn{2}{|c|}{Variable reference set} \\
& & & \multicolumn{2}{|c|}{reference set $m=3$} & \multicolumn{2}{|c|}{}\\
\cline{4-7}
$x_f$ & $M$  & $\varepsilon$ & Error & $N_f$ & Error & $N_f$\\
 \hline\hline
$2\pi$ & 10 & $10^{-5}$ & 0.0247309 & 300 & 6.47998e-05 & 550 \\
$2\pi$ & 10 & $10^{-9}$ & 0.0246415 & 480 & 2.24345e-09 & 1050\\
$4\pi$ & 10 & $10^{-5}$ & 0.888217 & 534 & 0.000142862 & 966\\
$4\pi$ & 20 & $10^{-9}$ & 0.0496889 & 960 & 1.05491e-08 & 2100\\
$6\pi$ & 10 & $10^{-5}$ & 14.4197 & 762 & 6.23799e-05 & 1530\\
$6\pi$ & 40 & $10^{-9}$ & 0.0232977 & 1560  & 3.06542e-09 & 3640\\
 \hline\hline
\end{tabular}
\end{center} 

\vspace*{-0.5cm}
\caption{Results for Example 2.}\label{tab2}
\end{table}

Now we compare the results obtained using equidistant nodes and Chebyshev points of second kind
for the reference set. For the same example the obtained results are given in Table \ref{tab3}.

\begin{table} [h]
\begin{center}
\begin{tabular}{|c|c|c|c|c|c|c|}
\hline\hline
& & & \multicolumn{2}{|c|}{Fixed equidistant}  & \multicolumn{2}{|c|}{Chebyshev fixed } \\
& & & \multicolumn{4}{|c|}{reference set $m=5$} \\
\cline{4-7}
$x_f$ & $M$  & $\varepsilon$ & Error & $N_f$ & Error & $N_f$\\
 \hline\hline
$2\pi$ & 10 & $10^{-5}$ & 6.93002e-05 & 400 & 2.69646e-05 & 400 \\
$2\pi$ & 10 & $10^{-9}$ & 1.91509e-05 & 650 & 8.13527e-06 & 650\\
$4\pi$ & 10 & $10^{-5}$ & 0.00215349 & 600 & 0.000338729 & 551\\
$4\pi$ & 20 & $10^{-9}$ & 3.85763e-05 & 1300 & 1.6391e-05 & 1300\\
$6\pi$ & 10 & $10^{-5}$ & 0.0275954 & 900 & 0.0164587 & 820\\
$6\pi$ & 40 & $10^{-9}$ & 1.00764e-05 & 2200  & 4.18516e-06 & 2200\\
 \hline\hline
\end{tabular}
\end{center} 

\vspace*{-0.5cm}
\caption{Results for Example 2.}\label{tab3}
\end{table}

As expected, the results using Chebyshev points of second kind are better than 
that obtained using equidistant nodes, due to the better approximation property of Lagrange interpolation polynomial with 
Chebyshev points of second kind toward the equidistant points, \cite{11}.

\item
 (\cite{3}, p. 245)
Keeping the differential system as in the previous example but changing the initial value conditions to $y(0)=[0.4,0,0,2],$ 
for $x_f=2\pi,\ M=20$ and $\varepsilon=10^{-9}$ with the method based on variable reference set
we obtained $\max_{0\le i\le M}\|y(x_i)-u_i\|=2.94126\cdot 10^{-9}$ and $N_f=1400.$ 

In this case the solution is $y(x)=[\cos{u}-0.6, \frac{-\sin{u}}{1-0.6 \cos{u}}, 0.8\sin{u}, \frac{0.8\cos{u}}{1-0.6\cos{u}}],$
where $x=u-0.6\sin{u}.$ 

Based on the previous examples the method with variable number of reference points is more efficient than the method with fixed number
reference points,
but we cannot deduce such a conclusion from the given convergence results.
\end{enumerate}

Using the method for stiff problems presented above we solved:

\begin{enumerate}
\setcounter{enumi}{3}
\item
$$
\left\{\begin{array}{lclcl}
\dot{y_1} & = & 998y_1+1998 y_2,&\quad& y_1(0)=1, \qquad x\in[0,1],\\
\dot{y_2} & = & -999y_1-1999 y_2,&\quad& y_2(0)=0, \\
\end{array}\right.
$$
with the solution $y_1=2e^{-x}-e^{-1000x}, y_2=-e^{-x}+e^{-1000x}.$

For $\tau=10$ the results are given in Table \ref{tab4}.

\begin{table} [h]
\begin{center}
\begin{tabular}{|c|c|c|c|c|c|}
\hline\hline
& & \multicolumn{2}{|c|}{Fixed equidistant}  & \multicolumn{2}{|c|}{Chebyshev fixed } \\
& & \multicolumn{4}{|c|}{reference set $m=5$} \\
\cline{3-6}
$M$  & $\varepsilon$ & Error & $N_f$ & Error & $N_f$\\
 \hline\hline
300 & $10^{-5}$ & 0.00164977 & 8585 & 0.000402419 & 8435 \\
500 & $10^{-7}$ & 0.000128781 & 10700 & 4.35037e-05 & 10555\\
 \hline\hline
\end{tabular}
\end{center} 

\vspace*{-0.5cm}
\caption{Results for Example 4.}\label{tab4}
\end{table}

\item

$$
\begin{array}{lcl}
\dot{y} &= & -20y, \qquad y(0)=1,\quad x\in[0,1].
\end{array}
$$

For $\tau=10, M=20$ and $\varepsilon=10^{-7}$ the results are given in Table \ref{tab5}.

\begin{table} [h]
\begin{center}
\begin{tabular}{|c|c|c|c|}
\hline\hline
 \multicolumn{2}{|c|}{Fixed equidistant}  & \multicolumn{2}{|c|}{Chebyshev fixed } \\
\multicolumn{4}{|c|}{reference set $m=5$} \\
\cline{1-4}
 Error & $N_f$ & Error & $N_f$\\
 \hline\hline
 1.19382e-06 & 800 & 4.58431e-07 & 785 \\
 \hline\hline
\end{tabular}
\end{center} 

\vspace*{-0.5cm}
\caption{Results for Example 5.}\label{tab5}
\end{table}
\end{enumerate}

To make the results reproducible we provide some code 
at \url{https://github.com/e-scheiber/scilab-ivpsolvers.git}.

\end{document}